\documentclass[12pt]{article}

\usepackage{graphicx,tikz,color,url}
\usepackage{amsmath,amssymb,latexsym}
\usepackage{pgfplots}
\usepackage{mathrsfs}
\usetikzlibrary{arrows}

\newtheorem{theorem}{Theorem}[section]

\newtheorem{lemma}[theorem]{Lemma}
\newtheorem{proposition}[theorem]{Proposition}

\newtheorem{corollary}[theorem]{Corollary}

\newcommand{\proof}{\noindent{\bf Proof.\ }}
\newcommand{\qed}{\hfill $\square$ \bigskip}
\newcommand{\cp}{\,\square\,}
\newcommand{\gp}{{\rm gp}}

\newcommand{\sr}{{}_{\rm SR}}

\newcommand{\ttt}{{}_{\rm tt}^{-}}
\newcommand{\diam}{{\rm diam}}
\newcommand{\gpt}{\gp_{\rm t}}
\newcommand{\gpd}{\gp_{\rm d}}
\newcommand{\gpo}{\gp_{\rm o}}

\begin{document}

\title{General position problems in strong and lexicographic products of graphs}

\author{Pakanun Dokyeesun$^{a,}$\thanks{\texttt{papakanun@gmail.com}}
\and Sandi Klav\v zar $^{a,b,c,}$\thanks{\texttt{sandi.klavzar@fmf.uni-lj.si}}
    \and Dorota Kuziak $^{d,}$\thanks{\texttt{dorota.kuziak@uca.es}}
\and Jing Tian $^{e,a,}$\thanks{\texttt{jingtian526@126.com}}
}
\maketitle

\begin{center}
$^a$ Institute of Mathematics, Physics and Mechanics, Ljubljana, Slovenia\\
\medskip

$^b$ Faculty of Mathematics and Physics, University of Ljubljana, Slovenia\\
\medskip

$^c$ Faculty of Natural Sciences and Mathematics, University of Maribor, Slovenia\\
\medskip

$^d$ Departamento de Estad\'istica e Investigaci\'on Operativa, Universidad de C\'adiz, Algeciras, Spain \\
\medskip
	
$^e$ School of Science, Zhejiang University of Science and Technology, Hangzhou, Zhejiang 310023, PR China\\
\medskip

\end{center}

\begin{abstract}
Outer, dual, and total general position sets are studied on strong and lexicographic products of graphs. Sharp lower and upper bounds are proved for the outer and the dual general position number of strong products and several exact values are obtained. For the lexicographic product, the outer general position number is determined in all the cases, and the dual general position number in many cases. The total general position number is determined for both products. Along the way some results on outer general position sets are also derived. 
\end{abstract}


\noindent
{\bf AMS Subj.\ Class.\ (2020)}:  05C12, 05C69, 05C76

\section{Introduction}

Let $G = (V (G), E(G))$ be a graph and $X \subseteq V (G)$. Vertices $u, v \in V (G)$ are {\em $X$-positionable} if for any shortest $u,v$-path $P$ we have $V(P)\cap X = \{u,v\}$. Clearly, adjacent vertices are $X$-positionable no matter whether they belong to $X$. Setting $\overline{X} = V (G)\setminus X$, the set $X$ is
\begin{itemize}
\item a {\em general position set}, if every $u, v \in X$ are $X$-positionable;
\item a {\em total general position set}, if every $u, v \in V(G)$ are $X$-positionable;
\item an {\em outer general position set}, if every $u, v \in X$ are $X$-positionable, and every $u \in X$, $v \in \overline{X}$ are $X$-positionable; and
\item  a {\em dual general position set}, if every $u, v \in X$ are $X$-positionable, and every $u, v \in \overline{X}$ are $X$-positionable.
\end{itemize} 

The cardinality of a largest general position set, a largest total general position set, a largest outer general position set, and a largest dual general position set of a graph $G$ will be respectively
denoted by $\gp(G)$, $\gpt(G)$, $\gpo(G)$, and $\gpd(G)$. These graph invariants will be respectively called the {\em general position number}, the {\em total general position number}, the {\em outer general position number}, and the {\em dual general position number} of $G$. Moreover, for any  $\tau(G) \in \{\gp(G), \gpt(G), \gpo(G), \gpd(G)\}$, by a {\em $\tau$-set} we mean a set of cardinality $\tau(G)$ having the respective property.

The above described variety of general position sets was introduced in~\cite{tian-2024+}, while (standard) general position sets have already been investigated a lot after their independent introduction in~\cite{chandran-2016, manuel-2018a}. In particular, a lot of attention has been given to general position sets of graph operations~\cite{ghorbani-2021, irsic-2023, klavzar-2021, korze-2023, tian-2021a, tian-2021b}. Among other contributions on the topic we highlight the papers~\cite{patkos-2020, thomas-2024, yao-2022}. With respect to graph operations, general position sets were also investigated in~\cite{klavzar-2019} on strong products. Among other matters covered in~\cite{klavzar-2019} we point to the problem whether the general position number of the strong product of two graphs equals the product of the general position numbers of the factors, which is one of the main outstanding problems in the area.  

In this paper we focus on the recently introduced three general position numbers of strong and  lexicographic products. In the next section we give necessary definitions and call up the known results that we need later. In Section~\ref{sec:more-outer} we prove some results on outer general position sets that are on one hand interesting in their own right and are on the other hand applied later on. The strong product is the topic of Section~\ref{sec:strong}. We prove general lower and upper bounds on the outer and the dual general position number, demonstrate their sharpness, and give exact values in some specific cases. Then, in Section~\ref{sec:lexicographic}, we turn our attention to lexicographic products. We first determine the outer general position number in all the cases, and complete the paper by establishing the dual general position number in many cases. Along the way we also determine the total general position number of both products.

\section{Preliminaries}

Unless stated otherwise, graphs $G=(V(G),E(G))$ considered in this paper are simple and connected. For a positive integer $k$, we use $[k]$ to represent the set $\{1,\ldots, k\}$. 
For a vertex $u$ of $G$, $N_G(u)$ denotes the set of neighbors of $u$ in $G$ and $N_G[u] = N_G(u)\cup \{u\}$. Vertices $u$ and $v$ of $G$ are {\em true twins} if $N_G[u] = N_G[v]$. True twins are thus adjacent vertices such that for each vertex $w\ne u,v$ we have $uw\in E(G)$ if and only if $vw\in E(G)$. 
The {\em degree} of a vertex $u$ is $d_G(u) = |N_G(u)|$ and the largest degree in the graph $G$ is the {\em maximum degree} $\Delta(G)$. A vertex of $G$ is {\em universal} if it is adjacent to all the other vertices of $G$.  
A vertex $u$ of $G$ is {\em simplicial} if $N_G(u)$ induces a complete subgraph. The set of simplicial vertices of $G$ will be denoted by $S(G)$ and the cardinality of $S(G)$ by $s(G)$. The order, the number of leaves, and the complement of $G$ are respectively denoted by $n(G)$, $n_1(G)$, and $\overline{G}$. Moreover, $\omega(G)$ and $\alpha(G)$ stand for the clique number and the independence number of $G$. 

The {\em distance} $d_G(u,v)$ between vertices $u$ and $v$ of $G$ is the number of edges on a shortest $u,v$-path.  The {\em diameter} of $G$ is the maximum distance between pairs of vertices of $G$ and is denoted by $\diam(G)$. Let $H$ be a subgraph of $G$. Then $H$ is {\em isometric} if for each pair of vertices $u,v\in V(H)$ we have $d_H(u,v) = d_G(u,v)$, and $H$ is {\em convex} if for any vertices $u,v\in V(H)$, any shortest $u,v$-path in $G$ lies completely in $H$. By abuse of language we will also say that a set $X\subseteq V(G)$ is convex if $X$ induces a convex subgraph. 

The vertex $u$ of $G$ is \emph{maximally distant} from a vertex $v\in V(G)$ if every $w\in N_G(u)$ satisfies $d_G(v,w)\le d_G(u,v)$. 
If $u$ is maximally distant from $v$, and $v$ is maximally distant from $u$, then $u$ and $v$ are \emph{mutually maximally distant}, MMD for short. Note that true twins are MMD. 
The set of MMD vertices of $G$ is called the {\em boundary} of $G$ and denoted by $\partial(G)$, see~\cite{bresar-2008, caceres-2005, kuziak-2018}, that is, 
$$\partial(G) = \{u\in V(G):\ u\  {\rm is\ an\ MMD\ vertex}\}\,.$$
We further set $b(G) = |\partial(G)|$. 
The {\em strong resolving graph} $G\sr$ of $G$ has $V(G)$ as the vertex set, two vertices being adjacent in $G\sr$ if they are MMD in $G$. The notion of the strong resolving graph was introduced in~\cite{oellermann-2007} while investigating the strong metric dimension of graphs. Note that there could be vertices in a graph $G$ which are not MMD with any other vertex in $G$ (they do not belong to $\partial(G)$). Hence such vertices are isolated in $G\sr$. In this sense, from now on we denote by $G\sr'$ the graph obtained from $G\sr$ by removing all (possibly zero) of its isolated vertices. Note that $V(G\sr')=\partial(G)$.

We next collect characterizations of the three newly introduced general position sets as respectively proved in~\cite[Theorems~2.1, 2.3, 3.1]{tian-2024+}.

\begin{theorem} \label{thm:all}
If $G$ is a connected graph and $X\subseteq V(G)$, then the following hold.
\begin{enumerate}
\item[(i)] $X$ is a total general position set of $G$ if and only if $X\subseteq S(G)$. Moreover, $\gpt(G)=s(G)$.
\item[(ii)] If $|X|\geq 2$, then $X$ is an outer general position set of $G$ if and only if each pair of vertices from $X$ is mutually maximally distant. Moreover, $\gpo(G) = \omega(G\sr)$.
\item[(iii)] If $X$ is a general position set of $G$, then $X$ is a dual general position set if and only if $G-X$ is convex.
\end{enumerate}
\end{theorem} 

Notice that by the definition of $G\sr'$ it follows that $\omega(G\sr')=\omega(G\sr)$. Moreover, by Theorem~\ref{thm:all}(ii), we have that for any connected graph $G$, 
\begin{align}
\gpo(G) = \omega(G\sr') \label{eq:3}\,. 
\end{align}

\section{On outer general position sets}
\label{sec:more-outer}

In this section we prove some results on outer general position sets that are interesting in their own right and that we will need later in the paper. The first of them holds for all types of general position sets. 

\begin{lemma}
\label{lem:isometric}
Let $H$ be an isometric subgraph of a graph $G$. If $X$ is (dual, outer, total) general position set of $G$, then $X\cap V(H)$ is a (dual, outer, total) general position set of $H$. 
\end{lemma}

\proof
Set $X_H = X\cap V(H)$. 

Assume that $X$ is a general position set of $G$ and suppose on the contrary that $X_H$ is not a general position set of $H$. Then there exist vertices $h,h',h''$ of $X_H$ such that $h''$ lies on a shortest $h,h'$-path $P$ in $H$. Since $H$ is isometric in $G$, the path $P$ is also a shortest path in $G$, but then $X$ is not a general position set of $G$. This contradiction shows that $X_H$ is a general position set of $H$. 

Assume now that $X$ is a dual general position set of $G$. We already know that $X_H$ is a general position set of $H$, hence in view of Theorem~\ref{thm:all}(iii) we claim that $V(H)\setminus X_H$ is convex in $H$. Assume on the contrary that this is not the case. Then there exist two vertices $u,v\in V(H)\setminus X_H$ and a shortest $u,v$-path $Q$ in $H$ which contains a vertex of $X_H$. As $H$ is isometric, $Q$ is also a shortest path in $G$. But since $u,v\in V(G)\setminus X$, this implies that $V(G)\setminus X$ is not convex in $G$, a contradiction with Theorem~\ref{thm:all}(iii). It follows that $X_H$ is a dual general position set of $H$. 

Assume next that $X$ is an outer general position set of $G$. Then $X_H$ is a general position set of $H$ and similarly as above we infer that also no shortest path between a vertex from $X_H$ and a vertex from $V(H)\setminus X_H$ contains a vertex of $X_H$. It follows that $X_H$ is an outer general position set of $H$. 

Finally, if $X$ is a total general position set of $G$, then the above arguments imply that $X_H$ is a total general position set of $H$. 
\qed

The case of Lemma~\ref{lem:isometric} for general position sets was earlier implicitly proved within the proof of~\cite[Theorem 3.1]{manuel-2018a}.

By $G\ttt$ we represent the graph obtained from $G$ by removing the edge between each pair of true twins. In particular, if $G$ has no true twins, then $G\ttt \cong G$, and $(K_n)\ttt$ is the edgeless graph of order $n$. 

\begin{proposition}
\label{prop:outer-diam-2}
If $G$ is a graph with $\diam(G) = 2$, then $\gpo(G) = \omega(G\ttt)$. In particular, if $G$ has no true twins, then $\gpo(G) = \alpha(G)$. 
\end{proposition}

\proof
We have already observed that every pair of true twins are MMD in a graph. 
Since $\diam(G) =2$, vertices $u$ and $v$ of $G$ are MMD if and only if either $d_G(u,v)=2$, or $u$ and $v$ are true twins. Then for every $u,v \in V(G)$, vertices $u,v$ are MMD in $G$ if and only if they are MMD in $G\ttt$ which means that $G\sr \cong (G\ttt)\sr$. By Theorem~\ref{thm:all}(ii) and $\diam(G) =2$, we have $\gpo(G) = \omega(G\sr) =\omega((G\ttt)\sr) = \omega(\overline{G\ttt}) =\alpha(G\ttt)$.

Note that any two non-adjacent vertices are MMD as $\diam(G) =2$. If $G$ has no true twins, we have $\omega(G\sr) = \omega(\overline{G}) = \alpha(G)$ and hence $\gpo(G) = \alpha(G)$.
\qed

A parallel result to Proposition~\ref{prop:outer-diam-2} for the case of general position sets was proved in~\cite[Proposition 2.4]{klavzar-2019}, where it is established that if $G$ has no true twins and $\diam(G) = 2$, then $\gp(G) = \omega(G\sr)$ if and only if $\gp(G) = \alpha(G)$.

If $k\ge 1$, then a set $X\subseteq V(G)$ is {\em $k$-independent} if $d_G(u,v) > k$ for each pair of vertices $u,v\in X$. The {\em $k$-independence number} of $G$ is the cardinality of a largest $k$-independent set in $G$ and denoted by $\alpha_k(G)$. Note that $\alpha_1(G) = \alpha(G)$.

\begin{proposition}
\label{prop:gpo-versus-alpka_k}
If $G$ is a graph with $\diam(G) = k \ge 2$, then $\gpo(G) \ge \alpha_{k-1}(G)$.
\end{proposition}

\proof
Let $X$ be a $(k-1)$-independent set of $G$ with $|X| = \alpha_{k-1}(G)$. Let $u,v\in X$, where $u\ne v$. As $\diam(G) = k$ and $d_G(u,v) > k-1$, we necessarily have $d_G(u,v) = k$. Using again the fact that $\diam(G) = k$ we see that $u$ and $v$ are MMD. It follows that $X$ induces a clique of $G\sr$, so that $\omega(G\sr) \ge |X| = \alpha_{k-1}(G)$. Theorem~\ref{thm:all}(ii) completes the argument. 
\qed

The bound of Proposition~\ref{prop:gpo-versus-alpka_k} is sharp for all $k\ge 2$. To see it, we give two constructions depending on the parity of $k$. 

First, let $k\ge 2$ be even and let $r = \frac{k-2}{2}$. Then the graph $K_{1,s}^r$, $s\ge 2$, is obtained from $K_{1,s}$ by subdividing all edges $r$ times. Notice that $\diam(K_{1,s}^r) = k = 2r+2$ and all leaves of this graph form the largest $(2r+1)$-independence set. Hence $\alpha_{k-1}(K_{1,s}^r) = n_1(K_{1,s}^r)$. Moreover, by Theorem~\ref{thm:all}(ii) we have that $\gpo(K_{1,s}^r) = \omega((K_{1,s}^r)\sr) = n_1(K_{1,s}^r)$. 

For $k\ge 3$ odd, let $t = \frac{k-1}{2}$. Then the graph $K_n^t$, $n\ge 2$, is obtained by attaching a path $P_t$ to every vertex of $K_n$. By construction, $\diam(K_n^t) = k = 2t+1$ and $\alpha_{k-1}(K_n^t) = n_1(K_n^t)$. Using Theorem~\ref{thm:all}(ii) again we can conclude that $\gpo(K_n^t) = \omega((K_n^t)\sr) = n_1(K_n^t) = 
\alpha_{k-1}(K_n^t)$.

\section{Strong products}
\label{sec:strong}

The \emph{strong product} of two graphs $G$ and $H$ is the graph $G \boxtimes H$ such that $V(G \boxtimes H) = V(G) \times V(H)$, and two vertices $(g,h), (g',h') \in V(G \boxtimes H)$ are adjacent if either $g=g'$ and $hh' \in E(H)$, or $gg' \in E(G)$ and $h=h'$, or $gg' \in E(G)$ and $hh' \in E(H)$. 
For a vertex $h\in V(H)$, set $G^h = \{(g,h)\in V(G\boxtimes H):\ g\in V(G)\}$. The set $G^h$ is called a {\em $G$-layer} of $G\boxtimes H$ and induces a subgraph of $G\boxtimes H$ isomorphic to $G$. For $g\in V(G)$, {\em $H$-layer} $^gH$ is defined as $^gH = \{(g,h)\in V(G\boxtimes H):\ h\in V(H)\}$. 
If $X\subseteq V (G\boxtimes H)$, the {\em projection} $p_G(X)$ of $X$ to $G$ is the set $\{g\in V(G):\ (g, h)\in X~{\rm for~some}~h\in V(H)\}$. The projection $p_H(X)$ of $X$ to $H$ is defined analogously. 

The following formulas are well-known, cf.~\cite{hammack-2011}: 
\begin{align}
d_{G \boxtimes H} ((g,h), (g',h')) & = \max\{d_G(g,g'), d_H(h,h')\} \quad  {\rm and} \label{eq:1} \\
N_{G \boxtimes H} [(g,h)] &= N_G[g] \times N_H[h]\label{eq:2}\,.  
\end{align}

In~\cite{klavzar-2019}, several exact results and bounds were proved on the general position number of strong products. In addition, the problem was posed whether $\gp(G\boxtimes H) = \gp(G) \gp(H)$ holds for any graphs $G$ and $H$~\cite[Problem 4.8]{klavzar-2019}. In this section we consider the other three versions of general position sets on strong products. 

The next result follows from~\eqref{eq:2} and can also be deduced from the proof of~\cite[Proposition 4]{caceras-2010}. 

\begin{lemma}
\label{lem:simplicial-strong}
If $G$ and $H$ are graphs, then $S(G\boxtimes H) = S(G)\times S(H)$. 
\end{lemma}

Lemma~\ref{lem:simplicial-strong} together with Theorem~\ref{thm:all}(i) yield: 

\begin{corollary}
If $G$ and $H$ are graphs, then $\gpt(G\boxtimes H) = s(G)s(H)$.
\end{corollary}

\subsection{Outer general position sets}
\label{sec:strong-outer}

To prove bounds on the outer general position number of strong products, we recall the following result that clarifies MMD vertices of strong products. 

\begin{lemma}{\rm \cite[Lemma 2.6]{oellermann-2007}}
\label{lem:mmd}  
Let $G$ and $H$ be two connected graphs. Let $g$, $g'$ be two vertices of $G$ and $h$, $h'$ two vertices of $H$. Then $(g,h)$ and $(g',h')$ are MMD in $G\boxtimes H$ if and only if one the following conditions holds:
\begin{enumerate}
\item[(i)] $g$, $g'$ are MMD in $G$ and $h, h'$ are MMD in $H$;
\item[(ii)] $g$, $g'$ are MMD in $G$ and $h = h'$;
\item[(iii)] $h$, $h'$ are MMD in $H$ and $g = g'$;
\item[(iv)] $g$, $g'$ are MMD in $G$ and $d_G(g,g') > d_H(h,h')$;
\item[(v)] $h$, $h'$ are MMD in $H$ and $d_G(g,g') < d_H(h,h')$. 
\end{enumerate}
\end{lemma}

Recall that the boundary $\partial(G)$ of $G$ contains its  MMD vertices and that $b(G) = |\partial(G)|$. The bounds on the outer general position number of strong products now read as follows. 

\begin{theorem}
\label{thm:strong}
If $G$ and $H$ are connected graphs of order at least $2$, then 
$$\gpo(G) \gpo(H)\le \gpo(G\boxtimes H) \le b(G) b(H)\,.$$ 
Moreover, if $G$ and $H$ are block graphs, then the bounds coincide. 
\end{theorem}

\proof
To prove the lower bound, let $X$ be a $\gpo$-set of $G$ and $Y$ be a $\gpo$-set of $H$. We will prove that $X\times Y$ is an outer general position set of $G\boxtimes H$. 
Consider any two vertices $(g,h), (g',h')$ from $X\times Y$. Assume first that $g\ne g'$ and $h\ne h'$. Since $X$ and $Y$ are $\gpo$-sets of $G$ and $H$, Theorem~\ref{thm:all}(ii) implies that $g, g'$ are MMD in $G$ and $h, h'$ are MMD in $H$. By Lemma~\ref{lem:mmd}(i), the vertices $(g,h), (g',h')$ are MMD in $G\boxtimes H$. Assume second that $g = g'$ and $h \ne h'$. Since $h, h' \in Y$ and $Y$ is a $\gpo$-set of $H$, the vertices $h, h'$ are MMD in $H$ by Theorem~\ref{thm:all}(ii). In view of Lemma~\ref{lem:mmd}(iii), the vertices $(g,h)$ and $(g',h')$ are MMD in $G\boxtimes H$. Analogously, if $h=h'$ and $g\ne g'$, then $(g,h), (g',h')$ are MMD in $G\boxtimes H$ by Lemma~\ref{lem:mmd}(ii). 
From Theorem~\ref{thm:all}(ii), we can conclude that $X \times Y$ is an outer general position set of $G\boxtimes H$, hence $\gpo(G\boxtimes H)\geq \gpo(G) \gpo(H)$.

To prove the upper bound, let $S$ be a $\gpo$-set of $G \boxtimes H$ and consider the following cases. 

\medskip\noindent 
{\bf Case 1}: $p_G(S)\subseteq \partial(G)$ and $p_H(S)\subseteq \partial(H)$. \\
In this case $\gpo(G\boxtimes H) = |S|\leq b(G) b(H)$. 

\medskip\noindent 
{\bf Case 2}: $S$ contains a vertex $(g,h)$ such that $h\notin \partial(H)$. \\
Let $(g',h')$ be a vertex from $S$, where $g'\ne g$. Since $S$ is an outer general position set of $G\boxtimes H$,  Theorem~\ref{thm:all}(ii) implies that $(g,h)$ and $(g',h')$ are MMD in $G\boxtimes H$. Since $h\notin \partial (H)$,  Lemma~\ref{lem:mmd} yields that $g$ and $g'$ are MMD in $G$. We thus see that each vertex of $p_G(S)\setminus \{g\}$ is adjacent to $g$ in the strong resolving graph $G\sr$ and hence  $|p_G(S)\setminus \{g\}|\leq \Delta(G\sr)$. Using Lemma~\ref{lem:mmd} once more we see that each $H$-layer $^{g'}H$ has at most $\gpo(H)$ vertices belonging to $S$, that is, $|V(^{g'}H)\cap S|\leq \gpo(H)$. Moreover, since $h\notin \partial(H)$, we have $|V(^{g}H)\cap S| = 1$. Indeed, for otherwise a vertex $(g,h')\in S$, where $h'\ne h$, would imply that $h\in \partial(H)$ by Lemma~\ref{lem:mmd}(iii). Then $|S|\leq \Delta(G\sr) \gpo(H) +1$.
By definitions of boundary sets and strong resolving graphs, we have $\Delta(G\sr) \le b(G)-1$ and $\omega(H\sr) \le b(H)$. Together with Theorem~\ref{thm:all}(ii) we get
\begin{align*}
   |S| & \leq \Delta(G\sr) \gpo(H) + 1 \\
   & \leq (b(G)-1) \omega(H\sr) +1 \\
   &\le (b(G)-1) b(H) +1 \\
   & < b(G) b(H). 
\end{align*}

\medskip\noindent 
{\bf Case 3}: $S$ contains a vertex $(g,h)$ such that $g\notin \partial(G)$. \\
By the commutativity of the strong product we can use parallel arguments as in Case 2 to derive the conclusion $|S|\leq \Delta(H\sr) \gpo(G) +1  < b(G) b(H)$.  

\medskip\noindent 
{\bf Case 4}: $S$ contains a vertex $(g,h)$ such that $g\notin \partial(G)$ and $h\notin \partial(H)$. \\
Lemma~\ref{lem:mmd} implies that one of $g$ and $h$ must actually respectively lie in $\partial(G)$ and $\partial(H)$, hence this case is not possible. 

\medskip
We have exhausted all the possibilities for $S$, hence the claimed upper bound is proved. 

\medskip
Assume now that $G$ and $H$ are block graphs. Then we know from~\cite{tian-2024+} that $\gpo(G) = s(G)$ and $\gpo(H) = s(H)$.  Since in a block graph every two simplicial vertices are MMD, we get $b(G) = s(G)$ and $b(H) = s(H)$. Hence the bounds of the theorem coincide if $G$ and $H$ are block graphs. 
\qed

The bounds of Theorem~\ref{thm:strong} are thus sharp on block graphs. On the other hand, recall from~\cite{sonnemann-1974} that $\alpha(C_5\boxtimes C_5) = 5$. By~\eqref{eq:1} we get $\diam(C_5\boxtimes C_5) = 2$, hence  Proposition~\ref{prop:outer-diam-2} yields $\gpo(C_5\boxtimes C_5) = 5$, the value which lies strictly between $\gpo(C_5)\gpo(C_5) = 4$ and $b(C_5)b(C_5) = 25$.  

Denoting the strong product of $k$ copies of $G$ by $G^{k,\boxtimes}$, the last result of this subsection reads as follows. 

\begin{proposition}
\label{prop:strong-trees}
If $G$ is a diameter two graph without true twins and $k\ge 1$, then $\gpo(G^{k,\boxtimes}) =  \alpha(G^{k,\boxtimes})$.
\end{proposition}

\proof
If $k=1$, the assertion is just Proposition~\ref{prop:outer-diam-2}. Since $G$ has no true twins, we can deduce by applying \cite[Lemma 9]{barragan-2016} that $G\boxtimes G$ also has no true twins. By induction on $k$, the power  $G^{k,\boxtimes}$ also has no true twins for each $k\ge 2$. Moreover, by \eqref{eq:1} and induction we also have $\diam(G^{k,\boxtimes}) = 2$. Applying Proposition~\ref{prop:outer-diam-2} again we reach the conclusion. 
\qed

\subsection{Dual general position sets}
\label{sec:strong-dual}

For the dual general position number of strong products we have the following general bounds.

\begin{theorem}
\label{thm:strong-dual-bounds}
If $G$ and $H$ are connected graphs, then 
\begin{align*}
s(G) s(H) \le \gpd(G\boxtimes H)\leq \min \{ & s(G)n(H) + s(H)n(G)-s(G)s(H),\\
& n(G)\gpd(H), n(H)\gpd(G)\}\,.
\end{align*}
\end{theorem}

\proof 
In~\cite[Corollary~3.2]{tian-2024+} it was proved that if $G$ is a graph and $X \subseteq S(G)$, then $X$ is a dual general position set of $G$. Hence $S(G\boxtimes H)$ is a dual general position set of $G\boxtimes H$. Lemma~\ref{lem:simplicial-strong} then yields the lower bound. 

Lemma~\ref{lem:isometric} applied to dual general position sets implies that $\gpd(G\boxtimes H)\leq n(G)\gpd(H)$ and $\gpd(G\boxtimes H)\leq n(H)\gpd(G)$. 

To prove that $\gpd(G\boxtimes H) \leq s(G)n(H) + s(H)n(G)-s(G)s(H)$, suppose on the contrary that there exists a dual general position set $R$ of $G \boxtimes H$ such that $|R| > s(G)n(H) + s(H)n(G)-s(G)s(H).$ Then there exists a vertex $(x,y)\in R$ such that $x\notin S(G)$ and $y\notin S(H)$, that is, $x$ and $y$ are not simplicial vertices of $G$ and $H$, respectively. Let $x', x'' \in N_G(x)$, and let $y', y'' \in N_H(y)$, where $x'x''\not\in E(G)$ and $y'y''\not\in E(H)$. Then the subgraph induced by $\{x',x,x''\} \times \{y',y,y''\}$ is isomorphic to $P_3 \boxtimes P_3$ with the universal vertex $(x,y)$. 

Since $R$ is a dual general position set of $G\boxtimes H$ and because $(x,y)\in R$, Theorem~\ref{thm:all}(iii) implies that one of the vertices $(x',y)$ and $(x'',y)$ must be in $R$, we may without loss of generality assume that $(x',y)\in R$. Using again that $R$ is a general position set of $G\boxtimes H$, we get $(x'',y')\notin R$ and $(x'',y'')\notin R$. Since the complement of $R$ is convex, we then infer that $(x',y')$ and $(x',y'')$ must be in $R$. But then the vertices $(x',y')$, $(x',y)$, and $(x',y'')$ induce a shortest $(x',y'),(x',y'')$-path, a contradiction with the assumption that $R$ is a general position set. 
\qed

To demonstrate that the three upper bounds of Theorem~\ref{thm:strong-dual-bounds} are pairwise incomparable, consider the following example. Let $C_n^+$ be the graph obtained from $C_n$ by attaching a new vertex to an arbitrary vertex of $C_n$. Now we consider $C_{2k+1}^+\boxtimes C_{2\ell+1}^+$, $k, \ell \ge 2$, and compute the values of the upper bounds of Theorem~\ref{thm:strong-dual-bounds}. Since  $s(C_{2k+1}^+) = s(C_{2\ell+1}^+) = 1$, the first bound is 
\begin{align*}
& \quad\ s(C_{2k+1}^+)n(C_{2\ell+1}^+) + s(C_{2\ell+1}^+)n(C_{2k+1}^+) - s(C_{2k+1}^+)s(C_{2\ell+1}^+) \\
& = (2\ell +2) + (2k+2) - 1 \\
& = 2\ell + 2k +3 \,.
\end{align*}
For the other two bounds, we first infer that $\gpd(C_{2k+1}^+) = 3$ and $\gpd(C_{2\ell+1}^+) = 3$. Assume without loss of generality that $k\leq \ell$. Then we have 
$$\min\{\gpd(C_{2k+1}^+)n(C_{2\ell+1}^+),\gpd(C_{2\ell+1}^+)n(C_{2k+1}^+)\} = \min \{6\ell+6, 6k+6\} = 6k +6\,.$$
Theorem~\ref{thm:strong-dual-bounds} thus gives
$$
\gpd(C_{2k+1}^+\boxtimes C_{2\ell+1}^+)\leq 
\begin{cases} 2\ell + 2k +3; & k \leq \ell\leq 2k+\frac{3}{2}, \\
6k +6; & \ell > 2k+\frac{3}{2}\,.
\end{cases}
$$

We conclude this section with the following exact result on the dual general position number of strong products. 

\begin{theorem}
\label{thm: strong-dual-complete}
If $H$ is a connected graph, then $\gpd(K_n\boxtimes H) = n \cdot \gpd(H)\,.$
\end{theorem}

\proof
By Theorem~\ref{thm:strong-dual-bounds}, $\gpd(K_n\boxtimes H) \le n\cdot \gpd(H)$. 

To prove the lower bound, let $V(K_n) = [n]$ and let $X$ be a $\gpd$-set of $H$. We claim that $X' = [n]\times X$ is a dual general position set of $K_n\boxtimes H$. We first need to show that $X'$ is a general position set of $K_n\boxtimes H$. For the sake of contradiction suppose this is not the case and let $(i,h)$, $(i',h')$, $(i'',h'')$ be three vertices of $X'$ such that $(i'',h'')$ lies on a shortest $(i,h),(i',h')$-path $P$ in $K_n\boxtimes H$. Then $p_H(P)$ induces a shortest $h,h'$-path in $H$ which contains $h''$, a contradiction with the assumption that $X$ is a general position set of $H$. 

We second need to verify that $V(K_n\boxtimes H)\setminus X'$ is convex in $K_n\boxtimes H$. For this sake consider arbitrary two vertices $(i,h)$ and $(i',h')$ from $V(K_n\boxtimes H)\setminus X'$. Suppose that there exists a shortest $(i,h),(i',h')$-path $P$ containing a vertex from $X'$, say $(j,k)$, such that $k\in X$. Then it follows that there exists a shortest $h,h'$-path containing the vertex $k$ in $H$. Since $V(H)\setminus X$ is convex, $d_H(h,h') < d_H(h,k) + d_H(k,h')$. By~\eqref{eq:1}, we have  
\begin{align*}
d_{K_n\boxtimes H}((i,h),(i',h')) & = d_H(h,h'),\\
d_{K_n\boxtimes H}((i,h),(j,k)) & = d_H(h,k),\\
d_{K_n\boxtimes H}((j,k),(i',h')) & = d_H(k,h').
\end{align*}
Then 
$$d_{K_n\boxtimes H}((i,h),(i',h')) < d_{K_n\boxtimes H}((i,h),(j,k)) + d_{K_n\boxtimes H}((j,k),(i',h'))\,,$$
which implies that $P$ is not a shortest $(i,h),(i',h')$-path in $K_n\boxtimes H$. This contradiction yields that $V(K_n\boxtimes H)\setminus X'$ is convex.
By Theorem~\ref{thm:all}(iii) we can conclude that $X'$ is a dual general position set of $K_n\boxtimes H$ and hence $\gpd(K_n\boxtimes H)\geq |X'| = n\cdot \gpd(H)$, and we are done.
\qed

Let $H$ be an arbitrary graph with $\gpd(H) = s(H)$, for instance, an arbitrary block graph. Then Theorem~\ref{thm: strong-dual-complete} yields 
$$\gpd(K_n\boxtimes H) = n \cdot \gpd(H) = s(K_n) s(H)\,,$$
which demonstrates sharpness of the lower bound of Theorem~\ref{thm:strong-dual-bounds}. 

\section{Lexicographic products}
\label{sec:lexicographic}

The \emph{lexicographic product} $G\circ H$ of two graphs $G$ and $H$ has $V(G \circ H) = V(G) \times V(H)$. Two vertices $(g,h), (g',h') \in V(G \circ H)$ are adjacent if either $g=g'$ and $hh' \in E(H)$, or $gg' \in E(G)$. The projections $p_G(X)$ and $p_H(X)$ of $X\subseteq V (G\circ H)$ on $G$ and $H$ are defined just as for the strong product, and $G$-layers and $H$-layers are also defined analogously. 

The following result is parallel to Lemma~\ref{lem:simplicial-strong}. Since we did not find it in the literature, its proof is provided. 

\begin{lemma}
\label{lem:simplicial-lex}
If $G$ and $H$ are connected graphs of order at least $2$, then 
$$
S(G\circ H) = 
\begin{cases} S(G)\times V(H); & H \text{ is complete}, \\
\emptyset; & \text{otherwise}\,.
\end{cases}
$$ 
\end{lemma}
\proof
We first assume that $H$ is a non-complete graph and suppose on the contrary that $S(G \circ H) \neq \emptyset$. Let $(x,y) \in S(G \circ H)$. As $G$ is a connected graph with at least two vertices, there exists $x' \in V(G)$ such that $xx' \in E(G)$. Since $H$ is non-complete, there exist two non-adjacent vertices $y', y'' \in V(H)$. Then $(x', y'), (x', y'') \in N_{G\circ H}((x,y))$ but they are not adjacent, a contradiction. 

Assume second that $H$ is complete. Consider a vertex $(x,y) \in S(G)\times V(H)$ and let $(x',y'), (x'',y'') \in N_{G\circ H}((x,y))$. If $x'=x''$, then $(x',y')$ and $(x'',y'')$ are adjacent since $H$ is complete. And if $x' \neq x''$, then $(x',y')$ and $(x'',y'')$ are adjacent because $x\in S(G)$ and hence $y'y''\in E(H)$. This proves that $S(G)\times V(H)\subseteq S(G\circ H)$. To see that also $S(G\circ H) \subseteq S(G)\times V(H)$ holds, we can use the argument from the first paragraph of the proof. 
\qed

Lemma~\ref{lem:simplicial-lex} together with Theorem~\ref{thm:all}(i) yield: 

\begin{corollary}
If $G$ and $H$ are graphs of order at least $2$, then
$$\gpt(G\circ H) =
\begin{cases} s(G)n(H); & H \text{ is complete}, \\
0; & \text{otherwise}\,.
\end{cases}
$$
\end{corollary}

\subsection{Outer general position sets}
\label{sec:lex-outer}

In this subsection we determine the outer general position number of lexicographic products. For this sake, we first need to recall some concepts and results presented in~\cite{kuziak-2016}. Moreover, the fact $\omega(G\circ H)=\omega(G)\omega(H)$, see~\cite[Theorem~3]{Doslic-2013}, is used several times in this subsection. 

Given a graph $G$, by $G_{\overline{2}}$ it is denoted the graph with vertex set $V(G_{\overline{2}})=V(G)$ such that two vertices $u,v$ are adjacent in $G_{\overline{2}}$ if either $d_G(u,v)\ge 2$ or $u,v$ are true twins. Further, by $G'$ we denote the graph obtained from $G$ by removing all its isolated vertices. (Note that this notation is consistent with $G\sr'$ used earlier for strong resolving graphs.) The next result is a direct consequence of~\cite[Remark 6]{kuziak-2016}, where $G+H$ denotes the join of graphs $G$ and $H$.

\begin{proposition}\label{prop:omega-G*}
 Let $G$ be a connected graph of order at least $2$.
\begin{enumerate}
\item[(i)] If $G$ has no universal vertex, then $\omega(G_{\overline{2}})=\omega((K_1+G)\sr')$.
\item[(ii)] If $\diam(G)\le 2$, then $\omega(G_{\overline{2}}')=\omega(G\sr')$.
\item[(iii)] If $G$ has no true twins, then $\omega(G_{\overline{2}})=\alpha(G)$.
\end{enumerate}
\end{proposition}

We next collect~\cite[Propositions 8, 11, 13, 16]{kuziak-2016} into the following result. Before that some extra preparation is needed. The {\em TF-boundary} of a non-complete graph $G$ is a set $\partial_{TF}(G) \subseteq \partial(G)$, where $x\in \partial_{TF}(G)$ whenever there exists $y\in \partial (G)$, such that $x$ and $y$ are MMD in $G$ and $N_G[x]\ne N_G[y]$ ($x,y$ are not true twins). Moreover, the \emph{strong resolving TF-graph} of $G$ is a graph $G_{SRS}$ with vertex set $V(G'_{SRS}) = \partial_{TF}(G)$, where two vertices $x,y$ are adjacent in $G'_{SRS}$ if $x$ and $y$ are MMD in $G$ and $N_G[x]\ne N_G[y]$. Finally, if $F$ is a subgraph of $H$, then for a given vertex $v\in V(G)$, by $(F)_v$ we denote the subgraph of $G\circ H$ isomorphic to $F$ which lies in the $H$-layer of $G\circ H$ corresponding to the vertex $v$ of $G$.  
 
\begin{proposition}{\rm \cite[Propositions 8, 11, 13, 16]{kuziak-2016}}
\label{prop:lex-SR}
If $G$ and $H$ are connected graphs of order at least $2$, then the following properties hold. 
\begin{enumerate}
\item[(i)] If $G$ has no true twins and $H$ is non-complete, then
$$(G\circ H)\sr'\cong \left(G\sr'\circ H_{\overline{2}}\right) \cup \bigcup_{v\in V(G)\setminus \partial(G)} (H_{\overline{2}}')_v. $$
\item[(ii)] If $m\ge 2$, then
$$(G\circ K_m)\sr'\cong (G\sr'\circ K_m) \cup \bigcup_{v\in V(G)\setminus \partial(G)} (K_m)_v.$$
\item[(iii)] If $m\ge 2$ and $H$ has no universal vertex, then
$$(K_m\circ H)\sr'\cong  \bigcup_{v\in V(K_m)} (H_{\overline{2}})_v.$$
\item[(iv)] If $G$ is non-complete and $H$ has no universal vertex, then
$$(G\circ H)\sr'\cong (G'_{SRS}\circ H_{\overline{2}}) \cup \bigcup_{v\in V(G)\setminus \partial_{TF}(G)} (H_{\overline{2}})_v.$$
\end{enumerate}
\end{proposition}

With Proposition~\ref{prop:lex-SR} in hand, we can deduce the subsequent results. In the first one, we avoid true twins in the first factor of the product, and complete graph in the second one.

\begin{theorem}
\label{thm:gpo-lex-H*}
Let $G$ be a connected graph without true twins of order at least $2$ and let $H$ be a non-complete graph.
\begin{enumerate}
\item[(i)] If $H$ has no universal vertex, then $\gpo(G\circ H)=\gpo(G) \gpo(K_1+H)$.
\item[(ii)] If $\diam(H) = 2$, then $\gpo(G\circ H)=\gpo(G) \gpo(H)$.
\item[(iii)] If $H$ has no true twins, then $\gpo(G\circ H)=\gpo(G) \alpha(H)$.
\end{enumerate}
\end{theorem}

\proof
By~\eqref{eq:3} and Proposition~\ref{prop:lex-SR}(i), we have 
\begin{align*}
\gpo(G\circ H) & =  \omega((G\circ H)\sr') = \omega(\left(G\sr'\circ H_{\overline{2}}\right) \cup \bigcup_{v\in V(G)\setminus \partial(G)} (H_{\overline{2}}')_v) \\
& = \omega\left(G\sr'\circ H_{\overline{2}}\right)=\omega(G\sr')\omega(H_{\overline{2}})\,.
\end{align*}
Combining Proposition~\ref{prop:omega-G*} with~\eqref{eq:3} we get the conclusion. 
\qed

\begin{theorem}
\label{thm:gpo-lex-H-complete}
If $G$ is a connected graph of order at least $2$ and $m\ge 2$, then
$$\gpo(G\circ K_m)=m \cdot \gpo(G).$$
\end{theorem}

\proof
By~\eqref{eq:3} and Proposition~\ref{prop:lex-SR}(ii), we have 
\begin{align*}
\gpo(G\circ K_m) & = \omega((G\circ K_m)\sr') = \omega(\left(G\sr'\circ K_m\right) \cup \bigcup_{v\in V(G)\setminus \partial(G)} (K_m)_v) \\ 
& = \omega\left(G\sr'\circ K_m\right) = \omega(G\sr')\omega(K_m).
\end{align*}
Using~\eqref{eq:3} once more, we get the formula. 
\qed

We have considered the case in which the second factor in the product is a complete graph. Since this product is not commutative, we now study the case in which the first factor is a complete graph.

\begin{theorem}
\label{thm:gpo-lex-G-complete}
Let $m\ge 2$ be an integer and let $H$ be a graph of order at least $2$ without universal vertex.
\begin{enumerate}
\item[(i)] If $\diam(H)=2$, then $\gpo(K_m\circ H)=\gpo(H)$.
\item[(ii)] If $\diam(H) > 2$, then $\gpo(K_m\circ H)=\gpo(K_1+H)$.
\end{enumerate}
\end{theorem}

\proof
By~\eqref{eq:3} and Proposition~\ref{prop:lex-SR}(iii), we have 
$$\gpo(K_m\circ H) = \omega((K_m\circ H)\sr') = \omega(\bigcup_{v\in V(K_m)} (H_{\overline{2}})_v) = \omega(H_{\overline{2}}).$$ Therefore, by Proposition~\ref{prop:omega-G*} and again~\eqref{eq:3} we conclude the proof.
\qed

Theorem~\ref{thm:gpo-lex-H*} deals with  first factors in the lexicographic product without true twins. Hence in the last result we focus on the case in which this first factor has true twins. 

\begin{theorem}
\label{thm:gpo-lex-SRS}
Let $G$ be a connected non-complete graph, and let $H$ be a graph of order at least $2$ without universal vertices.
\begin{enumerate}
\item[(i)] If $\diam(H)=2$, then $\gpo(G\circ H)=\omega(G'_{SRS})\gpo(H)$.
\item[(ii)] If $\diam(H) > 2$, then $\gpo(G\circ H)=\omega(G'_{SRS})\gpo(K_1+H)$.
\end{enumerate}
\end{theorem}

\proof
Using~\eqref{eq:3} and Proposition~\ref{prop:lex-SR}(iv), we get
\begin{align*}
\gpo(G\circ H) & = \omega((G\circ H)\sr') = \omega((G'_{SRS}\circ H_{\overline{2}}) \cup \bigcup_{v\in V(G)\setminus \partial_{TF}(G)} (H_{\overline{2}})_v) \\
& = \omega\left(G'_{SRS}\circ H_{\overline{2}}\right) = \omega(G'_{SRS})\omega(H_{\overline{2}})\,,
\end{align*}
and then Proposition~\ref{prop:omega-G*} and~\eqref{eq:3} yield the desired conclusion.
\qed

\subsection{Dual general position sets}
\label{sec:lex-dual}

We finally investigate the dual general position number. For this sake we recall that in~\cite[Theorem 2.1]{anand-2012}, non-complete convex sets of lexicographic products $G\circ H$ are characterized with three conditions. One of them asserts that $G\circ H$ can contain a non-complete convex subgraph only if $H$ is a complete graph. 
 
\begin{theorem}
If $G$ is a connected graph, then the following hold.
\begin{enumerate}
\item[(i)] If $G$ has no simplicial vertex and $H$ is a graph with no simplicial vertex, then $\gpd(G\circ H) = 0$.
\item[(ii)] If $n\ge 1$, then $\gpd(G\circ K_n) = \gpd(G)\cdot n$.
 
\end{enumerate}
\end{theorem}

\proof
(i) Let $X$ be a dual general position set of $G\circ H$ and let $(g,h)\in X$. Since $g\notin S(G)$, it has neighbors $g'$ and $g''$ such that $g'g''\notin E(G)$. Consider now the vertices $(g',h)$ and $(g'',h)$ and note that $d_{G\circ H}((g',h), (g'',h)) = 2$. Thus we cannot have $(g',h)\in X$ and $(g'',h)\in X$. Moreover, we also cannot have $(g',h)\notin X$ and $(g'',h)\notin X$ because $(G\circ H)\setminus X$ is convex. Hence, assume without loss of generality that $(g',h)\in X$ and $(g'',h)\notin X$. Let now $(g'',h')$ be a vertex with $d_H(h,h') = 2$. Such a vertex exists since $h$ is not simplicial. Then $(g',h)\in X$ and $(g,h)\in X$ imply that $(g'',h')\notin X$. But then $(g'',h')-(g,h)-(g'',h)$ is a shortest path with $(g,h)\in X$, a contradiction with the convexity of $(G\circ H)\setminus X$.

(ii) Let $X$ be a $\gpd$-set of $G$ and let $V(K_n) = [n]$. 
Set $Y = X\times [n]$ and $\overline{Y} = (V(G) - X)\times [n]$.
We first claim  that $Y$ is a general position set of $G\circ H$.
Suppose this is not the case and let $(g,h)$, $(g',h')$, $(g'',h'')$ be three vertices of $Y$ such that $(g'',h'')$ lies on a shortest $(g,h),(g',h')$-path $P$ in $G\circ K_n$. Then $p_G(P)$ induces a shortest $g,g'$-path in $G$ which contains $g''$, a contradiction with the assumption that $X$ is a general position set of $G$. 

Next, we prove that $\overline{Y}$ is convex. Suppose on the contrary that there exist two vertices $(g,h)$ and $(g',h')$ from $\overline{Y}$ such that some  shortest $(g,h),(g',h')$-path in $G\circ K_n$ contains a vertex $(g'', h'') \in Y$, where $g''\in X$. By the definition of $G\circ K_n$, since $\{h,h',h''\}$ induces a complete subgraph in $K_n$, our assumption implies that $g''$ lies on a shortest $g,g'$-path in $G$. By Theorem~\ref{thm:all}(iii), since $X$ is a dual general position set of $G$, 
there are no shortest $g,g'$-paths of $G$ containing the vertex $g''$. This contradiction implies that $\overline{Y}$ is convex. 

We have thus proved that $Y$ is a dual general position set of $G\circ K_n$ and hence $\gpd(G\circ K_n) \geq |Y| = \gpd(G)\cdot n$. 

Note that each $G$-layer has (at most) $\gpd(G)$ vertices contained in a dual general position set of $G\circ K_n$, hence we have $\gpd(G\circ K_n) \leq \gpd(G)\cdot  n$. We can conclude that $\gpd(G\circ K_n) = \gpd(G)\cdot n$.
\qed

\section*{Acknowledgements}

This work was supported by the Slovenian Research Agency ARIS (research core funding P1-0297 and projects N1-0285, N1-0218). D.\ Kuziak has been partially supported by the Spanish Ministry of Science and Innovation through the grant PID2023-146643NB-I00. D.\ Kuziak was visiting the University of Ljubljana supported by “Ministerio de Educaci\'on, Cultura y Deporte”, Spain, under the “Jos\'e Castillejo” program for young researchers (reference number: CAS22/00081).

\section*{Declaration of interests}

The authors declare that they have no known competing financial interests or personal relationships that could have appeared to influence the work reported in this paper.

\section*{Data availability}
Our manuscript has no associated data.


\end{document}